\documentclass{amsart}
\usepackage{amsfonts}
\usepackage{amssymb}
\numberwithin{equation}{section}
\usepackage{amscd}

\def\littlespacebox{\lower2ex\hbox{\vbox to 16pt{\vfill}}}
\def\spacebox{\lower2ex\hbox{\vbox to 22pt{\vfill}}}
\def\littlespacebox{\lower2ex\hbox{\vbox to 16pt{\vfill}}}

\newcommand{\A}{{\mathcal A}}
\newcommand{\B}{{\mathcal B}}
\newcommand{\F}{{\mathcal F}}

\newcommand{\CB}{{\mathcal C \mathcal B}}

\newcommand{\C}{{\mathbb C}}

\newcommand{\R}{{\mathbb R}}
\newcommand{\N}{{\mathbb N}}

\newcommand{\OL}{\mathcal {OL}}
\newcommand{\al}{\alpha}
\newcommand{\be}{\beta}
\newcommand{\de}{\delta}
\newcommand{\e}{\varepsilon}
\newcommand{\la}{\lambda}

\newcommand{\p}{\varphi}
\newcommand{\w}{\omega}
\newcommand{\ten}{\check\otimes}
\makeatletter
\mathchardef\hugecheck="7014
\newcommand\hugesize{\@setfontsize\hugesize{25pt}{0}}
\newcommand\smallhugesize{\@setfontsize\smallhugesize{20pt}{0}}

\def\specialchecksmall{\mbox{\hbox to
0pt{\raisebox{-4pt}{\smallhugesize$\hugecheck$}}
                                $\kern-2.7pt\otimes$}}

\makeatother

\newtheorem{claim}{}[section]
\newtheorem{theorem}[claim]{Theorem}
\newtheorem{lemma}[claim]{Lemma}
\newtheorem{proposition}[claim]{Proposition}

\newtheorem{remark}[claim]{Remark}

\newtheorem{example}[claim]{Example}

\def\proclaim #1. #2\par{\medbreak
\noindent{\bf#1.\enspace}{\sl#2}\par\medbreak}
\begin{document}

\title[Approximation Properties for  $L_p(VN(G))$ ]
{Approximation Properties for Non-commutative $L_p$-Spaces
Associated with Discrete Groups}
\date{April 15,  2001, Revised March 12, 2002}
\author[Marius Junge]{Marius Junge$^*$}
\address{Department of Mathematics\\
University of Illinois, Urbana, IL 61801 USA}
\email[Marius Junge]{junge@math.uiuc.edu}
\author[Zhong-Jin Ruan]{Zhong-Jin Ruan$^*$}
\address{Department of Mathematics\\
University of Illinois, Urbana, IL 61801 USA}
\email[Zhong-Jin Ruan]{ruan@math.uiuc.edu}
\thanks{${}^*$The authors were partially supported by the National Science
Foundation DMS-0088928 and DMS-9877157}

\begin{abstract}  Let  $1 < p < \infty$.
It is shown that if $G$ is a discrete group
with  the approximation property
introduced by Haagerup and Kraus,
then  the non-commutative  $L_p(VN(G))$ space has
the operator space approximation property.
If, in addition, the  group von Neumann algebra $VN(G)$
has the QWEP, i.e. is a quotient of a  $C^*$-algebra with Lance's 
weak expectation
property,  then
$L_p(VN(G))$ actually has the completely contractive approximation property
and the approximation maps can be chosen to be finite-rank completely 
contractive
multipliers on $L_p(VN(G))$.
Finally, we show that if  $G$ is a countable discrete group  having 
the approximation
property and  $VN(G)$ has the QWEP, then
$L_p(VN(G))$ has a very nice
local structure, i.e. it  is a $\mathcal C\OL_p$ space and has   a 
completely bounded
Schauder basis.
\end{abstract}

\maketitle
\let\text=\mbox
\let\cal=\mathcal

\section{Introduction}

Approximation properties for Banach spaces were first studied by
Grothendieck \cite {Gr}.
The corresponding non-commutative analogue of Grothendieck's program
has been developed and successfully  applied to operator algebras,
and   more recently to operator spaces.
These non-commutative approximation properties have played a crucial role in
the study of von Neumann algebras  and $C^*$-algebras (particularly,
for group von Neumann algebras and group $C^*$-algebras).
For example, it is well-known (from \cite {Lance}, \cite {Co} and \cite {EL})
that  a discrete group  $G$ is
amenable if and only if the reduced group $C^*$-algebra  $C^*_{red}(G)$
is  nuclear (respectively, the group von Neumann algebra $VN(G)$ is injective).
Some weaker conditions (i.e. weak amenability and approximation property)
for locally compact groups have been studied by Haagerup and Kraus 
(see \cite {Ha2}
and  \cite {HK}).
It was shown in \cite {HK} that a discrete group $G$ has the 
approximation property
(or simply, AP) if  and only if  $C^*_{red}(G)$ has the operator  space
approximation property (or simply, OAP) of Effros and Ruan \cite {ER1}
(respectively, $VN(G)$  has the weak$^*$ OAP
of Kraus  \cite {Kraus}).
The approximation property for a group $G$ is also closely related to
its Fourier algebra $A(G)$, which can be identified with the
operator predual of $VN(G)$.
Since   a dual operator space $V^*$ has the weak$^*$ OAP
if and only if $V$ has the OAP (see \cite {EKR}),
we can conclude that a discrete group $G$ has the AP if and only if its Fourier
algebra $A(G)$  has the OAP (see $\S 2$ for a direct proof).

Non-commutative $L_p$-spaces (for $1\le p < \infty$)  are also  very
important   in  the study of non-commutative harmonic analysis and
non-commutative   probability theory.
Given a von Neumann algebra $R$, we may isometrically identify  $L_1(R)$ with
the predual $R_*$ of $R$.
However, it is very important to note that since we usually use the
\emph{trace  duality} between $R$ and $R_*$, the correct
operator space matrix norm on
$L_1(R)$ should be given by the opposite operator space matrix norm on
$R_*$, i.e. we should completely isometrically identify $L_1(R)$ with 
$(R_*)^{op} =
(R^{op})_*$ (see detailed explanations in $\S 3$).
Then we may use Pisier's complex interpolation method to obtain a canonical
operator space matrix norm on the non-commutative $L_p(R)$ space.
The purpose of this paper is to study the approximation properties for
non-commutative $L_p(VN(G))$ spaces associated with discrete groups $G$.
Our main results can be stated as follows.

\begin{theorem}
Let $1< p< \infty$.
If $G$ is a discrete group with the AP, then  $L_p(VN(G))$
has the OAP.
\end{theorem}

Using a similar technique, we can also prove (in Proposition
\ref {P3.cb}) that
if $G$ is a weakly amenable discrete group and $1 < p < \infty$,
then $L_p(VN(G))$ has the completely bounded approximation property.

It is known by Grothendick \cite {Gr} (see \cite {LT}) that if  a
Banach space dual $V$ is separable  and has the Grothendieck's 
approximation property,
then $V$ actually  has the contractive approximation property.
The separability can be removed if $V$ is a reflexive space (see \cite {DU}).
At this moment, we can not obtain such a general result for
operator space duals.
However, we may obtain the following result with
an additional assumption that $VN(G)$ has the \emph{QWEP}
introduced by Kirchberg \cite {Ki1},
i.e. $VN(G)$ is a  quotient of a $C^*$-algebra with
Lance's  weak expectation property.
A $C^*$-algebra $\A$ is said to have  the \emph{weak expectation property}
if for the universal  representation $\pi: \A \hookrightarrow B(H)$,
there is a completely positive and contractive map $P: B(H) \to 
\A^{**}$ such that
$P \circ \pi =  id_\A$  (see Lance \cite {Lance}).

\begin{theorem}
Let $1 < p < \infty$.
If $G$ is a discrete group with the AP and  $VN(G)$ has the QWEP,
then $L_p(VN(G))$ has the completely contractive approximation property.
\end{theorem}

Combining Theorem 1.2 and the recent results in \cite {JNRX}, we can show
in Theorem \ref {P4.T} that  if $G$ is a countable discrete group
satisfying the conditions given in Theorem 1.2, then $L_p(VN(G))$ has 
a very nice
local structure, i.e. it is  a $\mathcal C \OL_p$ space and has a completely
bounded  Schauder basis (see definitions in $\S 5$).
This is a quite surprisng result, and it is only true for $1 < p < \infty$.
Indeed, for $p = 1$,  it is known from \cite {ER2} and \cite  {JNRX}
that  $L_1(VN(G))= VN(G)_*^{op}$ (equivalently, $VN(G)_*$) is a $\mathcal C
\OL_1$ space if and only if $VN(G)$ is an injective von Neumann algebra.
For $p = \infty$, we need to consider the reduced group $C^*$-algebra 
$C^*_{red}(G)$,
and it is known from \cite {ER2} that  $C^*_{red}(G)$ is a $\mathcal 
C \OL_\infty$ space
if and only if $C^*_{red}(G)$ is a nuclear $C^*$-algebra.
Therefore,  in the case of $p = 1$ or $\infty$, this can only happen when
a discrete group $G$ is  amenable.

The paper is organized  as follows.
In $\S 2$, we recall some necessary notions and  results  on
operator spaces and completely bounded multipliers of  Fourier algebras,
and we clarify some simultaneous convergence properties for completely bounded
multipliers of  $C^*_{red}(G)$ and $A(G)$  (see Propositions \ref{P2.2} and
\ref{P2.3}).
In $\S 3$, we recall the complex interpolation for operator spaces introduced
by Pisier \cite {PiOH}
and  recall non-commutative  $L_p(R)$ spaces arising from
von Neumann algebras $R$.
Since we are mainly interested in the case of  group von Neumann algebras
$ VN(G)$ for discrete groups $G$,  it suffices to consider von 
Neumann  algebras
with  normal faithful tracial states.
The readers are refered to \cite {Ko}, \cite {Te}, \cite {Iz} and \cite {Fi}
for non-commutative $L_p$-spaces arising, by complex interpolation,
from general von Neumann algebras.
We  prove  Theorem 1.1  in $\S 3$, and prove Theorem 1.2 in $\S 4$.
Motivated by an argument given in \cite [Theorem 2.1]{HK}, we
are able to show in Theorem \ref {P4.multiplier} that the completely 
contractive
approximation maps  in Theorem 1.2 can actually be chosen to be completely
contractive  multipliers on $L_p(VN(G))$.
This result will be very useful in non-commutative harmonic analysis.
We end the paper by $\S 5$, in which we study the local structure and
completely bounded  Schauder basis for $L_p(VN(G))$ spaces.
We include some interesting examples of countable residually finite
discrete groups $G$ with the AP, for which the $L_p(VN(G))$  spaces
are ${\mathcal C \OL}_p$ spaces and have completely bounded Schauder bases.

We wish to thank Quanhua Xu for discussions (with the second author)
on the subject and some details concerning the complex interpolation method
and approximation properties.

\section {Preliminaries}

We assume that the readers are familiar with the basic notions in
operator algebras and  operator spaces.
The readers are refered to
Takesaki's book \cite {Ta} and Stratila's book  \cite {St} for details
on operator algebras,
and are refered to  Paulsen's book \cite{Paulsen}, the recent book of
Effros and Ruan \cite{ERbook}, and Pisier's book \cite {Pibook}
for the details on operator spaces and completely bounded maps.
We recall that an operator space $V$ is said to have the \emph{completely
bounded approximation property} (or simply, \emph{CBAP})
if there exists a net of finite-rank maps
$T_\al : V \to V$ such that $\|T_\al\|_{cb} \le \lambda$ for some constant
$\lambda$ and  $T_\al \to id_V$ in the \emph{point-norm topology} on $V$, i.e.
we have $\|T_\al(x) - x\| \to 0$ for all $x \in V$.
We let
\begin{equation}
\label {F2.1}
\Lambda(V) = \inf \{\lambda\}
\end{equation}
denote the CBAP constant of $V$.
An operator space $V$ is said to have the \emph{ completely contractive
approximation property} (or simply, \emph{CCAP}) if $\Lambda(V) = 1$.
In \cite {ER1}, Effros and the second author studied the operator 
space analogue of
Grothendieck's  approximation property for Banach spaces.
We recall that an operator space $V$ is said to have the
\emph{operator space approximation property}
(or simply, \emph{OAP}) if  there exists a net of finite-rank maps
$T_\al : V \to V$ such that
$T_\al \to id_V$ in the \emph{stable point-norm topology} on $V$, i.e.
we have $\|T_\al\otimes id_\infty(x) - x\| \to 0$ for all $x \in 
V\ten K_\infty$,
where $\ten$ denotes the \emph{operator space injective tensor product}.
It is easy to see that CBAP implies OAP.

Let us assume that $G$ is a discrete group.
The \emph{left regular representation}  $\lambda : G\rightarrow B(\ell_{2}(G))$
is defined by
\begin{equation}
\label {F2.2}
\lambda (s)\xi (t)=\xi (s^{-1}t)
\end{equation}
for all $\xi \in \ell_{2}(G)$ and $s,t\in G$.
If we let $\lambda(\C [G]) = {\rm span}\{\la(s), ~ s\in G\}$, then the
\emph{reduced group $C^*$-algebra} $C^*_{red}(G)$ and the
\emph{group von Neumann algebra} $VN(G)$ are the norm closure and  the weak$^*$
closure of $\lambda(\C [G])$ in $B(\ell_{2}(G))$, respectively.
The \emph { Fourier algebra}
\[
A(G)  = \{f :  ~ f(t) = \langle \la (t) \xi ~ | ~ \eta\rangle ~
{\rm for ~ some}  ~ \xi, \eta \in \ell_2(G)\}
\]
is  the space of  all coefficient functions of the left regular representation
$\lambda $.
Given $f \in  A(G)$, its norm is given by
\[
\|f\|  =  \inf \{\|\xi\| \, \|\eta\|:
~ f (t)= \langle \la (t) \xi ~ | ~ \eta\rangle\}.
\]
It was shown by Eymard \cite {Ey}  that $A(G)$ with this norm and the pointwise
multiplication $m$ is a (contractive) commutative Banach algebra.
Moreover, there is a natural operator space matrix norm on $A(G)$ obtained by
identifying $A(G)$ with the operator predual $VN(G)_*$  of $VN(G)$.
With this operator space matrix norm, the multiplication
$m : A(G) \times A(G) \to A(G)$ is completely contractive
in the sense that
\[
\|[m(f_{ij}, h_{kl})]\| \le \|[f_{ij}]\|\, \|[h_{kl}]\|
\]
for all $[f_{ij}] \in M_m(A(G))$ and $[h_{kl}]\in M_n(A(G))$, and thus
$A(G)$ is a completely contractive Banach algebra (see \cite{ERCov}
and \cite {Ruan2}).

A function $\p$ on $G$ is  called a \emph{multiplier}
of $A(G)$ if  $m_\p (f) =  \p f$ maps $A(G)$ into $A(G)$.
It is easy to show from the closed graph theorem that if $\p$ is a 
multiplier of
$A(G)$, then the map $m_\p : A(G) \to A(G)$ is automatically  bounded.
If the map $m_\p$ is completely bounded on $A(G)$, we call $\p$
a \emph { completely bounded multiplier}  of $A(G)$.
We let $M_0A(G)$ denote the space of all completely bounded multipliers
of $A(G)$, which is  equipped with the cb-norm on $A(G)$.
It is clear that if $\p \in A(G)$, then the multiplication map
$m_\p : A(G) \to A(G)$  is completely  bounded with
\[
\|m_\p\|_{cb} \le \|\p\|.
\]
Therefore, we have a norm decreasing inclusion $A(G) \hookrightarrow M_0A(G)$.
If $G$ is an amenable group then  $A(G) \hookrightarrow M_0A(G)$ is
an isometric inclusion.

Amenability is one of the most important subjects in  harmonic analysis.
We recall that a group $G$ is amenable
if and only if  $A(G)$ has a contractive approximate identity.
If  we let $A_c(G)$ denote the space of all elements in $A(G)$ with
compact supports,  then $A_c(G)$ is norm dense in $A(G)$ and thus
the amenability of $G$ is equivalent to the existence of
a net of $\p_\al \in A_c(G)$ such that
$\|\p_\al \| \le 1$ and $m_{\p_\al}\to id_{A(G)}$
in the point-norm topology.
Haagerup introduced a weaker amenability condition for $G$ in  \cite {Ha2}.
Let us recall that a group $G$ is said to be  \emph {weakly amenable}
if there exists a net of  $\p_\al \in A_c(G)$ such that
$\|m_{\p_\al}\|_{cb} \le \la $  for some constant $\la $ and
$m_{\p_\al} \to id_{A(G)}$ in the point-norm topology.
It is clear from the definition that  amenability implies  weak amenability,
but the converse is not true (see  \cite {Ha1}, \cite {CH} and \cite {DH}).

It is also  known   that $M_0A(G)$ equipped with the cb-norm is a dual
space. It has a  predual $Q(G)$, which is the closure of
$\ell_1(G)$ under the norm given by
\[
\|f\|_{Q} = \sup\{|\int_G f(t) \p (t) \, dt | :
~ \p \in M_0A(G), \|\p\|_{cb} \le 1\}
\]
(see  \cite{He}, \cite{DH} and \cite[$\S 1$]{HK}).
A  group $G$ is said to have the
\emph{approximation property}
(or simply, \emph{AP}) if there exists a net of $\{\p_\al\} $ in $A_c(G)$
such that $m_{\p_\al} \to id _{A(G)}$ in the $\sigma(M_0A(G), Q(G))$ topology.
It was shown by Haagerup and Kraus \cite {HK} that weak amenability implies AP,
but the converse is not true (see $G = {\Bbb Z}^2\rtimes SL(2, {\Bbb Z})$).

Given $\p \in M_0A(G)$, the adjoint  map $M_\p = m_\p^* : VN(G) \to VN(G)$
is a weak$^*$ continuous  completely  bounded map  on $VN(G)$, which satisfies
\[
M_\p(\lambda(t) ) = \p(t) \lambda(t).
\]
Since $M_\p$ maps $C^*_{red}(G)$ into $C^*_{red}(G)$, it also induces a
completely bounded map $\overline M_\p: C^*_{red}(G) \to C^*_{red}(G)$.
All these maps have the same cb-norm, i.e. we have
\[
\|m_\p\|_{cb} =  \|M_\p\|_{cb} = \|\overline M_{\p}\|_{cb}.
\]
Let us write $M_\infty = B(\ell_2)$,  $K_\infty=K(\ell_2)$ and
$T_\infty=K(\ell_2)^* = B(\ell_2)_*$, and let $\overline \otimes$ and
$\hat \otimes$ denote the \emph{ normal spatial} and \emph{operator space
projective} tensor products, respectively.
Then for any $\p\in M_0A(G)$, we may obtain  completely bounded maps
$m_\p \otimes id_{\infty}$, $\overline M_\p\otimes id_\infty$ and
$M_\p \otimes id_\infty$ on $A(G)\ten {K_\infty}$,
$C^*_{red(G)} \ten {K_\infty}$ and $VN(G) \overline \otimes 
M_\infty$, respectively.
Given  $a \in VN(G)\overline \otimes M_\infty$ and $f \in A(G)\hat 
\otimes T_\infty$,
we may define a  bounded linear functional    $\w_{a, f}$ on $M_0A(G)$ given by
\begin{equation}
\label {F2.m1}
\w_{a, f}(\p) = \langle M_\p\otimes id_\infty(a), f\rangle.
\end{equation}
We may similarly define  bounded linear functionals on $M_0A(G)$ by considering
\begin{equation}
\label {F2.m2}
\tilde \w_{a, f}(\p) = \langle \overline M_\p\otimes id_\infty (a), \, f\rangle
\end{equation}
for $a \in C^*_{red}(G)\ten K_\infty$ and $f
\in (C^*_{red}(G)\ten K_\infty)^*$, and
\begin{equation}
\label {F2.m3}
\tilde{\tilde \w}_{a, f} (\p) = \langle a, ~   m_\p\otimes id_\infty 
(f) \rangle
\end{equation}
for $f \in A(G) \ten K_\infty$ and $a \in (A(G) \ten K_\infty)^* $.
Haagerup and Kraus proved in \cite [Propositions 1.4  and 1.5]{HK}
that the bounded linear functionals defined
in (\ref {F2.m1}) and (\ref {F2.m2}) are  all contained  in $Q(G)$,
and on the other hand, all linear functionals in $Q(G)$ have  such forms.
We note that this is also true for  (\ref {F2.m3}).
We summarize these in the following proposition.

\begin{proposition} {\rm (Haagerup and Kraus \cite {HK})}
\label {P2.1}
Let $G$ be a discrete group.  Then the bounded linear functionals defined in
{\rm (\ref {F2.m1})}, {\rm (\ref {F2.m2})} and {\rm (\ref {F2.m3})} are all
contained in $Q(G)$, and we actually have
\begin{eqnarray*}
Q(G) &=& \{\w_{a, f}:  ~ a \in VN(G)\overline \otimes M_\infty,
~ f \in A(G) \hat \otimes T_\infty\} \\
&=& \{\tilde \w_{a, f}:  a \in C^*_{red}(G)\ten K_\infty, ~ f
\in (C^*_{red}(G)\ten K_\infty)^*\} \\
&=& \{\tilde {\tilde \w}_{a, f}: f \in A(G) \ten K_\infty, ~ a \in
(A(G) \ten K_\infty)^*
\}.
\end{eqnarray*}
\end{proposition}
\begin{proof}
We  only need to study (\ref {F2.m3}) and the last equality.
Given  $f \in A(G) \ten K_\infty$,
$a \in (A(G) \ten K_\infty)^* = VN(G) \hat \otimes T_\infty$ and $\e > 0$,
it is known from the operator space theory (see \cite {ER1} and \cite{ERbook})
that  we may  write
\begin {equation}
\label {F2.4}
a = \al \tilde a \be = [\sum_{j, k} \alpha_{ij} \tilde a_{jk} \beta_{kl}]
\end{equation}
for some  Hilbert-Schmidt matrices $\al=[\al_{ij}], \be = [\be_{kl}] 
\in HS_\infty$ with
$\|\al\|_2 = \|\be\|_2 = 1$ and
$\tilde a = [\tilde a_{jk}] \in VN(G) \overline  \otimes M_\infty$
with $\|\tilde a \| < \|a\| + \e$.
Then for any $\p \in M_0A(G)$, we have
\begin{eqnarray*}
\tilde {\tilde \w}_{a, f}(\p)  &=& \langle a, ~ m_\p\otimes id_\infty 
(f) \rangle
= \langle \al \tilde a \be, ~ m_\p\otimes id_\infty (f)\rangle \\
&=&  \langle M_\p\otimes id_\infty (\tilde a) , \al^{tr} f \be^{tr} \rangle
= \w_{\tilde a, \al^{tr} f \be^{tr}}(\p).
\end{eqnarray*}
Since $\al^{tr} f \be^{tr} \in A(G) \hat \otimes T_\infty$, we can 
conclude from
the first equality that
$\tilde {\tilde \w}_{a, f} = \w_{\tilde a, \al^{tr}f\be ^{tr}}$ is an
element in $Q(G)$.

On the other hand, given any $\w \in Q(G)$, it is known from
\cite[Proposition 1.5]{HK}   that  we can write
$\w=  \w_{a, f}$ for  some  $a \in VN(G) \overline \otimes M_\infty$ and
$f \in A(G) \hat \otimes T_\infty$.
Using a similar caluculation as that given in (\ref {F2.4}), we may write
$f = \al \tilde f \be$  for some $\al, \be \in HS_\infty$ and
$\tilde f \in A(G) \ten K_\infty$, and thus
\[
\w = \w_{a, f} = \tilde {\tilde \w}_{\tilde a, \tilde f}
\]
with $\tilde a = \al^{tr} a \be ^{tr} \in VN(G) \hat \otimes T_\infty 
= (A(G)\check
\otimes K_\infty)^*$.
\end{proof}

  Proposition \ref {P2.1} shows that a net of completely bounded
multipliers $\{\p_\al\}$  converges to $\p$ in the  $\sigma(M_0A(G), 
Q(G))$ topology
if and only if the corresponding net of completely bounded maps 
$\{m_{\p_\al}\} $
  converges to $m_\p$
(respectively, $\{\overline M_{\p_\al}\} $ converges to $\overline M_{\p}$)
in the stable point-weak topology.
If $G$ is a discrete group, then every  element $\p$ in $A_c(G)$ has a
finite support,  and thus   $m_\p$ is a finite-rank map on $A(G)$
(respectively,  $\overline M_{\p}$ is a finite-rank map on $C^*_{red}(G)$).
It follows that $G$ has the AP if and only if there exists a net of
$\{\p_\al\}$ in $A_c(G)$ such that the finite-rank maps $m_{\p_\alpha}$
  converge to $id_{A(G)}$
(respectively,  $\overline M_{\p_\al}$ converge to
$id_{C^*_{red}(G)}$) in the stable point-weak topology.
By a standard convexity argument we may obtain a net of $\{\tilde \p_\al\}$
in $A_c(G)$ such that $m_{\tilde \p_\al} \to id_{A(G)}$
in the stable point-norm topology.
This gives a direct proof that $G$ has the AP if and only if $A(G)$ 
has the OAP.
Similarly, we can also find a (possibly different) net of $\{\tilde \psi_\be\}$
in $A_c(G)$ such that  $\overline M_{\tilde \psi_\be}\to
id_{C^*_{red}(G)}$ in the stable point-norm topology.
The following proposition shows that we can actually choose the same net
of $\{\p_\al\} $ in $A_c(G)$ such that $m_{\p_\al}\to id_{A(G)}$ and
$\overline M_{\p_\al}\to id _{C^*_{red}(G)}$ simultaneously.

\begin{proposition}
\label{P2.2}
Let $G$ be a discrete group with the AP.
Then there exists a net of $\{\p_\al\}$ in $A_c(G)$ such that
$m_{\p_\al} \to id_{A(G)}$ in the stable point-norm  topology on $A(G)$ and
$\overline M_{\p_\al} \to id_{C^*_{red}(G)}$ in the
stable point-norm  topology on $C^*_{red}(G)$.
\end{proposition}
\begin{proof}
It suffices to show that given any $f \in A(G) \ten K_\infty$,
$a \in C^*_r(G) \ten K_\infty$ and $\e > 0$, we can find an element
$\p\in A_c(G)$ such that
\[
\|m_{\p} \otimes id_\infty(f) - f\| < \e
~\, ~  {\rm  and }   ~\, ~
\|\overline M_{\p} \otimes id_\infty(a) - a\| < \e.
\]

If $F= \{f_1, \cdots, f_n\}$ is  any (non-empty) finite subset
of $(C^*_r(G) \ten K_\infty)^*$, then  we have $\tilde \w_{a, f_i}   \in Q(G)$,
and it follows from Proposition \ref {P2.1} that
there exist $\tilde a_i \in VN(G) \hat \otimes T_\infty$
and $\tilde f_i \in A(G) \ten K_\infty$ such that
\[
\tilde \w_{a, f_i}  = \tilde {\tilde \w}_{\tilde a_i, \tilde f_i}.
\]
Since $G$ has the AP,  there exists  $\p_F \in A_c(G)$ such that
\[
\|m_{\p_F} \otimes id_\infty(f) - f\| < \e
\]
and
\[
\|m_{\p_F} \otimes id_\infty(\tilde f_i) - \tilde f_i\| < \frac 1 {n 
(K+1)}  ~ \,
~ \, ~\, ({\rm for} ~ 1 \le i \le n),
\]
where we let $K = {\rm max}\{\|\tilde a_i\|\}$.
It follows that
\begin{eqnarray*}
|\langle \overline M_{\p_F} \otimes id_\infty(a) - a, ~ \, f_i\rangle |
&=& |\tilde \w_{a, f_i}(\p_F - 1)| = |\tilde {\tilde \w}_{\tilde a, 
\tilde f_i}(\p_F - 1)|\\
&= &|\langle \tilde a_i, ~ \, m_{\p_F} \otimes id_\infty(\tilde f_i) 
- \tilde f_i\rangle |
< \frac 1 {n}.
\end{eqnarray*}
Then we get a net of elements $\{\p_F\}$ in $A_c(G)$, which is
indexed by finite subsets $F$ of  $(C^*_r(G) \ten K_\infty)^*$, such that
\[
\|m_{\p_F} \otimes id_\infty(f) - f\| < \e
\]
for all $\p_F$ and
$\overline M_{\p_F} \to id_{C^*_{red}(G)}$ in the stable point-weak topology.
By a standard convexity argument, we may find an element  $\p\in 
A_c(G)$ such that
\[
\|m_{\p} \otimes id_\infty(f) - f\| < \e
~\, ~ {\rm  and}
~\, ~
\|\overline M_{\p} \otimes id_\infty(a) - a\| < \e.
\]
\end{proof}

Haagerup \cite {Ha2} proved that if $G$ is weakly amenable, then
\begin{equation}
\label {F2.5}
\Lambda(G) = \Lambda(A(G)) = \Lambda (C^*_{red}(G)) < \infty.
\end{equation}
Using a similar argument to that given  in Proposition \ref {P2.2}, we can
obtain  the following proposition for weakly amenable groups.

\begin{proposition}
\label {P2.3}
Let $G$ be a weakly amenable discrete group.
Then there exists a net of $\{\p_\al\}$
in $A_c(G)$ such that
\[
\|m_{\p_\al}\|_{cb} = \|\overline{M}_{\p_\al}\|_{cb} \le \Lambda(G)
\]
and $m_{\p_\al}\to id_{A(G)}$ and $\overline{M}_{\p_\al}\to 
id_{C^*_{red}(G)}$ in
the point-norm topologies on $A(G)$ and $C^*_{red}(G)$, respectively.
\end{proposition}

\begin{proof}
Let us outline the proof.
First, let us fix arbitrary  $f \in A(G)$ and $a \in C^*_{red}(G)$.
If $F = \{f_1, \cdots, f_n\}$ is a finite subset of $C^*_{red}(G)^*$,
we have $\tilde \w_{a, f_i} \in Q(G)$ and thus there exists $\tilde a_i \in
VN(G) \hat \otimes T_\infty$ and $\tilde f_i \in A(G) \check \otimes K_\infty$
such that
\[
\tilde \w_{a, f_i} = \tilde {\tilde \w}_{\tilde a_i, \tilde f_i}.
\]
Since $G$ is weakly amenable,  for any  $\e > 0$
there exists $\p_F \in A_c(G)$ such that
$\|m_{\p_F}\|_{cb} \le \Lambda(G) < \infty$ and
\[
\|m_{\p_F} (f) - f\| < \e
\]
and
\[
\|m_{\p_F} \otimes id_\infty (\tilde f_i) - \tilde f_i\| < \frac 1 {n 
(K+1)}  ~ \,
~ \, ~\, ({\rm for} ~ 1 \le i \le n),
\]
where we let $K = {\rm max}\{\|\tilde a_i\|\}$.
Then we can show as the proof given in Proposition \ref {P2.2}
that $\overline {M}_{\p_F}\to id _{C^*_{red}(G)}$ in the point-weak
topology, and thus there exists $\p$ in the convex hull of
$\{\p_F\}$ for which
\[
\|m_{\p}\|_{cb} = \|\overline {M}_{\p}\|_{cb} \le \Lambda(G)
\]
and
\[
\|m_{\p}(f) - f\| < \e ~\, ~  {\rm and} ~\, ~  \|\overline 
{M}_{\p}(a) - a\| < \e.
\]
A similar argument can be applied to arbitrary finite collection of elemets in
$A(G)$ and $C^*_{red}(G)$, and this completes the proof.
\end{proof}

\section{ Approximation Properties for $L_p(VN(G))$}

Let us first briefly recall some basic notions from the complex interpolation
theory (see \cite{BL}).
A pair of Banach  spaces $(V, W)$ is said to be a \emph{compatible couple}
if there is a  topological vector space $X$ and continuous inclusions
$V\hookrightarrow X$ and $W\hookrightarrow X$, which allow us to identify
$V \cap W$ and $V + W$ in $X$.
Then there is a canonical Banach space  norm on $V \cap W$ given by
\[
\|x\|_{V \cap W} = {\rm max } \{\|x\|_V, \|x\|_W\},
\]
and a canonical Banach space norm  on  $V+W$ given by
\[
\|x\|_{V+W} = \inf\{ \|v\|_{V} + \|w\|_{W}: ~ x = v+w \}.
\]

Let $S $ denote  the strip $\{z\in \C: ~ 0\le {\rm Re} z \le 1 \}$
in the complex plane $\C$ and let $S^0$ denote the interior of $S$.
We let $\F $ denote the collection of all continuous and bounded functions
$f : S \to V+W$ such that
\begin{itemize}
\item [(1)] $f$ is analytic on $S^0$,
\item [(2)] $f(it) \in V$ and $f(1+it) \in W$ for all  $t \in \R$,
\item[(3)] $f(it)\to 0$ and  $f(1+it) \to 0$ as  $t\to \infty$.
\end{itemize}
It is clear that $\F$ is a vector space.
Actually, $\F$ is a Banach space with the norm given
by
\[
\|f\|_{\F} = {\rm max}\{\sup\{\|f(it)\|_{V}:  ~ t \in \R\},~ 
\sup\{\|f(1+it)\|_{W}:
~ t \in \R\}\}.
\]
For $0\le \theta \le 1$, the space
\[
(V, W)_{\theta} = \{ a \in V+W: ~ a = f(\theta)  ~ {\rm for ~ some}~ f \in \F\}
\]
is called  the \emph{complex interpolation} of $(V, W)$.
This is a Banach space with the norm given by
\[
\|a\|_\theta = \inf\{\|f\|_\F:  ~ a = f(\theta), f \in \F \}.
\]
If $V$ and $W$ are operator spaces, Pisier \cite{PiOH} showed that 
there is a canonical
operator space matrix norm on $(V, W)_\theta$ given by
\begin{equation}
\label {F3.1}
M_n ((V, W)_\theta ) = (M_n(V), M_n(W))_\theta
\end{equation}
for all $n \in \N$.

Let $R$ be a von Neumann algebra with a normal faithful tracial state $\tau$.
For $1\le p <\infty$, the non-commutative  $L_p(R)$ space is defined
to be the closure of $R$ under the $p$-norm
\[
\|x\|_p = \tau((x^* x)^{\frac p2})^{\frac 1p}.
\]
We usually write  $L_\infty (R) = R$.
The trace $\tau$ induces  a canonical contractive embedding $j_1 : R \to R_*$
given by
\begin{equation}
\label {F3.4}
\langle j_1(x), ~  y\rangle =  \tau (xy).
\end{equation}
With this embedding, $(R, R_*)$ is a compatible couple of Banach spaces,
and we actually have the isometry
\begin{equation}
\label{F3.5a}
L_p(R) = (R, R_*)_{\frac 1p}
\end{equation}
for   $1 \le p < \infty$ (see Kosaki \cite{Ko}).

However, we have to be very careful about their operator space matrix norms.
In the operator space theory  it is custumary to use
the following \emph{parallel duality}  paring
\begin{equation}
\label {F3.e1}
\prec [\beta_{ij}], [\alpha_{ij}]\succ
= \sum_{i,j = 1}^n  \beta_{ij} \alpha_{ij} = Tr(\beta \alpha^{op})
\end{equation}
between $T_n = M_n^*$ and $M_n$.
The reason to use this parallel duality paring instead of the
\emph{trace duality} paring
\begin{equation}
\label {F3.e2}
\langle [\beta_{ij}], [\alpha_{ij}]\rangle
= \sum_{i,j = 1}^n \beta_{ij} \alpha_{ji} = Tr(\beta \alpha)
\end{equation}
is that we will be able to get the complete isometry
\[
M_n(V)^* = T_n \hat \otimes V^*
\]
for every operator space $V$.
Therefore, if we wish to use the trace duality paring (\ref  {F3.4}),
we should define
\[
L_1(R)  = (R^{op})_* = (R_*)^{op},
\]
for which we can obtain the complete isometry
\[
L_1(M_n\overline{\otimes} R) = T_n \hat \otimes L_1(R).
\]
The opposite operator space $(R_*)^{op}$ can be isometrically
identified with  $R_*$,
  but is equipped with the opposite operator space matrix norm
\[
\|[f_{ij}^{op}]\| = \|[f_{ji}]\|
\]
for all $[f_{ij}^{op}] \in M_n((R_*)^{op})$.
Then we may use (\ref {F3.1}) to  obtain a
canonical operator space   matrix norm  on $L_p(R)$ given by
\begin{equation}
\label {F3.e6}
M_n(L_p(R)) = (M_n(R), M_n(L_1(R)))_{\frac 1p}.
\end{equation}
The following proposition shows that we can replace $R$ in (\ref {F3.5a}) and
(\ref {F3.e6}) by a weak$^*$ dense  $C^*$-subalgebra of $R$.

\begin{proposition}
\label {P3.1}
Let $R$ be a von Neumann algebra with a normal faithful tracial state $\tau$
and let $\B$ be a weak$^*$ dense $C^*$-subalgebra of $R$.
For $1 < p < \infty$, we have the complete isometry
\[
L_p(R) = (\B, L_1(R))_{\frac 1p}.
\]
\end{proposition}
\begin{proof}
Let us first consider the isometry case.
Given any  $1 < p < \infty$, there exists a positive number $q$
such that $1 < q < p$.
It follows from the reiteration theorem (see \cite {BL} $\S 4.6$) that
we have the isometry
\[
L_p(R)= ((R, L_1(R))_0, (R, L_1(R))_{\frac 1q})_{\frac qp} = (R, 
L_q(R))_{\frac qp}.
\]
Since $L_q(R)$ is a reflexive space, we also have the isometry
\[
(\B, L_q(R))_{\frac qp}^* =(\B^*, L_q(R)^*)_{\frac qp}
= (L_{q'}(R), \B^*)_{1-\frac qp}
\]
(see \cite {BL} $\S 4.5$),
where $q'= {\frac q {q-1}}$ is the conjugate exponent to $q$.
Since $\B \subseteq R$, there is a canonical  contractive inclusion
\[
(\B, L_q(R))_{\frac qp} \to (R, L_q(R))_{\frac qp} = L_p(R),
\]
from which we obtain a canonical contraction (by taking the adjoint)
\[
L_{p'}(R) = L_{p}(R)^*\to  (L_{q'}(R), \B^*)_{1-\frac qp} .
\]

If we let $\F_0$ denote the space of elements having the form
\[
f(z) = exp (\la z^2) \sum _{n=1}^N exp (\lambda_n z) x_n
\]
with $\la > 0$, $N \in \Bbb N$, $\lambda_n \in \R $ and
$x_n \in L_{q'}(R) \cap \B^* \subseteq L_1(R)$,
then it is known from \cite{BL} that $\F_0$ is dense in $\F$
(with respect to the compatible couple $(L_{q'}(R), \B^*)$)
and thus the space of
elements $f(1-{\frac qp})$ with $f \in \F_0$ is norm dense in
$(L_{q'}(R), \B^*)_{1-\frac qp}$.
If $f \in \F_0$ then for every  $t\in \R$, we have
\begin{eqnarray*}
\|f(1+it)\|_{\B^*} &=& \sup \{|\tau(f(1+it)y)|: ~  y\in \B, \|y\| \le 1 \} \\
&=& \sup \{|\tau(f(1+it)y)|: ~ y \in R, \|y\| \le 1 \}
= \|f(it)\|_{L_1(R)}
\end{eqnarray*}
by the Kaplansky's density theorem.
This shows that  we actually have the isometry
\[
L_{p'}(R) = (L_{q'}(R), \B^*)_{1-\frac qp},
\]
and thus the isometry
\[
L_p(R) = (\B, L_q(R))_{\frac qp}.
\]

It was shown in Wolff \cite [Theorem 2]{Wolff}  that if we are given 
Banach spaces
$V_i$ ($i=1,2,3,4$) such that $V_1 \cap V_4$ is norm dense in $V_2$ 
and $V_3$, and
$(V_2, V_4)_\theta = V_3$, $(V_1, V_3)_\phi = V_2$ for some $0 < 
\theta, \phi < 1$,
then we have $(V_1, V_4)_\xi = V_2$ for
$\xi = \frac {\phi \theta } {1 - \phi + \phi \theta}$.
Now if we let $V_1 = \B, V_2 = L_p(R), V_3 = L_q(R)$ and $V_4 = 
L_1(R)$, and let
$\theta = \frac {p-q} {q(p-1)}$ and $\phi = \frac qp$, then we obtain
$\xi = \frac 1p$ and the isometry
\[
L_p(R) = (\B, L_1(R))_{\frac 1p}
\]
by  Wolff's result.

For each $n \in \Bbb N$, we have the contractive linear isomorphisms
\[
M_n((\B, L_1(R))_{\frac 1p}) = (M_n(\B), M_n(L_1(R)))_{\frac 1p}
\to (M_n(R), M_n(L_1(R)))_{\frac 1p}=M_n(L_p(R)),
\]
and
\[
M_n((\B, L_1(R))_{\frac 1 {p'}}) = (M_n(\B), M_n(L_1(R)))_{\frac 1{p'}}
\to (M_n(R), M_n(L_1(R)))_{\frac 1{p'}}=M_n(L_{p'}(R)).
\]
Then by duality,  we must have the isometry
\[
M_n(L_p(R)) = (M_n(\B), M_n(L_1(R)))_{\frac 1p}.
\]
\end{proof}

If $G$ is a discrete group, then we have a canonical normal faithful tracial
state  $\tau$ on $VN(G)$ given by
\begin{equation}
\label {F3.5}
\tau (x) = \langle x \de_e~  | ~ \de_e \rangle,
\end{equation}
where $\de_e$ is the characteristic function at the
unit element $e \in G$.
In this case, $(C^*_{red}(G), L_1(VN(G)))$ is a compatible
couple with the canonical embedding
$j_1: C^*_{red}(G) \to L_1(VN(G))$  given by
\begin{equation}
\label {F3.j1}
j_1 (\lambda(s)) = \delta_{s^{-1}} = \check \delta_s
\end{equation}
for all $s \in G$.
Given a function $\p : G \to \C$, we let $\check \p$ be the function on
$G$  defined by $\check \p(t) = \p(t^{-1})$.
For $1 < p < \infty$, we have the complete isometry
\[
L_p(VN(G)) = (C^*_{red}(G), L_1(VN(G)))_{\frac 1p}
\]
by Proposition \ref {P3.1}, and
have a  canonical  (contractive) inclusion
\[
j_p: C^*_{red}(G) \to L_p(VN(G)),
\]
for which the range space  $j_p(C^*_{red}(G))$ is  norm dense
in $L_p(VN(G))$.
If we let $p'$ be the conjugate exponent to $p$, then the duality
between  $L_p(VN(G))$ and $L_{p'}(VN(G))$ is given by
\begin{equation}
\label {F3.conj}
\langle j_p(\lambda(s)), ~ j_{p'}(\lambda(t))\rangle = \tau (\la(s) \la(t))
\end{equation}
for all $s, t \in G$.

In this case, the space $L_1(VN(G))$  can be explicitly expressed as follows.
First let us recall that there is a normal $*$-anti-automorphism  $\kappa$
on the group von Neumann algebra $VN(G)$ defined by
\[
\kappa(\sum_{i} a_i \lambda(t_i)) =  \sum _i a_i \lambda(t_i^{-1})
\]
(see  \cite {St} and \cite {ES}).
Since  $\kappa$ is a normal $*$-anti-automorphism on  $VN(G)$,
we must have
\[
\|\kappa_n ([x_{ij}])\| = \|[x_{ji}]\|
\]
for all $[x_{ij}] \in M_n(VN(G))$.
Its pre-adjoint map $\kappa_*$ induces an isometry on $A(G)$
which   satisfies
\[
\|(\kappa_*)_n ([f_{ij}])\| = \|[f_{ji}]\|
\]
for all $[f_{ij}] \in M_n(A(G))$.
Then we may completely isometrically identify  $L_1(VN(G))$ with
$\kappa_*(A(G)) = \{\kappa_*(f) = \check f:  ~ f \in A(G)\}$ and thus
obtain
\[
\|[\check f_{ij}]\| = \|[f_{ji}]\|
\]
for $[\check f_{ij}] \in M_n(L_1(VN(G)))$.

If  $\p$ is a completely bounded multiplier
on $A(G)$, we can prove that $\check \p$ is a completely bounded multiplier on
$L_1(VN(G))$ such that  $\|m_{\check \p}\|_{cb} = \|m_{\p}\|_{cb}$.
To see this, let us assume that $\p \in M_0A(G)$.
Then for any $\check f \in L_1(VN(G))$, we have
\[
m_{\check \p} ({\check f}) = {\check \p}{\check f} = \check {(\p f)}
= \kappa_*(\p f) \in L_1(VN(G)).
\]
This shows that
\begin{equation}
\label {F3.k}
m_{\check \p}= \kappa_*\circ m_\p \circ \kappa_*
\end{equation}
is a  well-defined multiplier map on $L_1(VN(G))$. Since
\begin{eqnarray*}
\|[m_{\check \p} ({\check f}_{ij})]\| & =& \|[\kappa_*(\p f_{ij})]\|
= \|[{\p f_{ji}}]\|  = \|[m_{\p} ({f}_{ji})]\|
\end{eqnarray*}
for  every $[\check f_{ij}] \in M_n(L_1(VN(G)))$,  we can conlcude that
\begin{equation}
\label {F3.n}
\|m_{\check \p}\|_{cb} = \|m_{\p}\|_{cb}.
\end{equation}
In general if $T$ is a completely bounded map (respectively, a 
complete isometry)
on $A(G)$, then $T$ induces a completely bounded
map (respectively, a complete isometry)
\[
\check T = \kappa_* \circ T \circ \kappa_*
\]
  on $L_1(VN(G))$, which satisfies
\[
\|\check T\|_{cb} = \|T\|_{cb}.
\]
As a consequence   we may replace the approximation maps $m_{\p_\alpha}$ on
$A(G)$ in  Proposition \ref {P2.2} and Proposition \ref {P2.3} by
the corresponding  maps  $m_{\check \p_\alpha}$ on  $L_1(VN(G))$
and obtain the following modified result.

\begin{proposition}
\label{P2.4}
If $G$ is a discrete group with the AP,
then there exists a net of $\{\p_\al\}$ in $A_c(G)$ such that
$m_{{\check \p}_\al} \to id_{L_1(VN(G))}$ in the stable point-norm
topology on $L_1(VN(G))$ and
$\overline M_{\p_\al} \to id_{C^*_{red}(G)}$ in the
stable point-norm  topology on $C^*_{red}(G)$.

If $G$ is a weakly amenable discrete group, then there exists a net
of $\{\p_\al\}$  in $A_c(G)$ such that
\[
\|m_{{\check \p}_\al}\|_{cb} = \|\overline{M}_{\p_\al}\|_{cb} \le \Lambda(G)
\]
and $m_{{\check \p}_\al}\to id_{L_1(VN(G))}$ and $\overline{M}_{\p_\al}\to
id_{C^*_{red}(G)}$ in the point-norm topologies on $L_1(VN(G))$ and
$C^*_{red}(G)$, respectively.
\end{proposition}

\begin{proposition}
\label {P3.2}
Given $\p \in M_0A(G)$, $(\overline M_{\p}, m_{\check  \p})$ is a
compatible pair of completely bounded maps on the compatible couple
$(C^*_{red}(G),  L_1(VN(G)))$.
\end{proposition}
\begin{proof}
Assume that we are given $\p \in  M_0A(G)$.
For any $s, t\in G$, we have
\[
\langle j_1(\overline M_{\p} (\lambda (s))), ~ \lambda (t) \rangle
= \langle j_1(\p(s) \lambda (s)), ~ \lambda (t) \rangle
= \p(s)  \tau (\lambda (s) \lambda (t))
\]
and
\[
\langle m_{\check  \p} (j_1(\la(s)), ~ \la(t)\rangle  = \p(t^{-1}) 
\tau(\la(s) \la(t))
= \p(s) \tau(\la(s) \la(t)).
\]
This shows that
\[
j_1(\overline M_{\p} (\lambda (s))) =  m_{\check  \p} (j_1(\la(s))
\]
for all $s \in G$.
Since $\la(\C [G])$ is norm dense in $C^*_{red}(G)$,
we obtain
\begin{equation}
\label {F3.7}
j_1\circ \overline M_{\p} =  m_{\check  \p}\circ j_1,
\end{equation}
and thus $(\overline M_{\p}, m_{\check  \p})$ is compatible pair  on
$(C^*_{red}(G),  L_1(VN(G)))$.
\end{proof}

\begin{remark}
\label {P3.4}
{\rm Let $\otimes ^h$ denote the Haagerup tensor product for operator spaces
(see \cite {ERbook}).
If $(V_i, W_i)$ with $(i = 1,2)$  are compatible couples of operator spaces,
Pisier proved in \cite{PiOH} Theorem 2.3 that
$(V_1 \otimes ^h V_2, W_1 \otimes ^h W_2)$
is again a compatible couple of operator spaces and we have the 
complete isometry
\begin{equation}
\label {F3.P1}
(V_1 \otimes ^h V_2, W_1 \otimes ^h W_2)_\theta = (V_1, W_1)_\theta \otimes ^h
(V_2, W_2)_\theta.
\end{equation}
It is also known that for arbitrary operator space $V$, we have the 
complete isometry
\begin{equation}
\label {F3.P2}
V \ten K_\infty = C\otimes ^h V \otimes ^h R,
\end{equation}
where $C$ and $R$ are \emph{column} and \emph{row} operator Hilbert spaces over
$\ell_2$.
It follows from (\ref {F3.P1}) and (\ref {F3.P2}) that  we have  the 
complete isometry }
\[
(V, W)_\theta \ten K_\infty = (V\ten K_\infty, W \ten K_\infty)_\theta
\]
{\rm  for arbitrary compatible couple of operator spaces $(V, W)$.
In particular for $1 < p < \infty$, we have the complete isometry}
\[
L_p(VN(G)) \ten K_\infty
= (C^*_{red}(G) \ten K_\infty, L_1(VN(G)) \ten K_\infty)_{\frac 1p}.
\]
\end{remark}

\noindent \emph{Proof of Theorem 1.1.}
We need to show that for every
$a \in L_p(VN(G)) \ten K_\infty$ and $\e > 0$,
there exists a finite-rank map $T$ on $L_p(VN(G))$ such that
\[
\| T \otimes id_\infty (a) - a\| < \e.
\]
Since  $L_p(VN(G)) \ten K_\infty
= (C^*_{red}(G) \ten K_\infty,  L_1(VN(G))  \ten K_\infty)_{\frac 1p}$,
there exists a continuous and bounded map
\[
f : S \to C^*_{red}(G) \ten K_\infty +  L_1(VN(G))  \ten K_\infty
\]
in  $\F$ such that $a = f({\frac 1p})$.
Since $f(it) \in C^*_{red}(G) \ten K_\infty$ and $f(it) \to 0$ as
$t\to \infty$,
the set $\{f(it)\}_{t \in \R}$ is contained in a
compact subset of $C^*_{red}(G) \ten K_\infty$.
Then there exists an element $(x_n) \in C^*_{red}(G) \ten K_\infty 
\ten c_0$ such that
$\{f(it)\}_{t\in \R}$ is contained in the norm closure of the convex hull of
$\{x_n\}$ in $C^*_{red}(G) \ten K_\infty$.
Similarly, there exists an element $(y_n) \in  L_1(VN(G))  \ten 
K_\infty \ten c_0$ such
that  the set $\{f(1+it)\}_{t\in \R}$ is contained in the norm 
closure of the convex hull
of
$\{y_n\}$ in $ L_1(VN(G))  \ten K_\infty$.
Since $G$ has the AP and we may identify $K_\infty \ten c_0$ with a 
closed subspace
of $K_\infty \ten K_\infty \cong K_\infty$, it follows from
Proposition \ref {P2.4} that  there exists $\p\in A_c(G)$ such that
\[
\|\overline M_{\p} \otimes id_\infty (f(it)) - f(it)\| < \e
\]
and
\[
\|m_{\check  \p} \otimes id_\infty (f(1+it)) - f(1+it)\| < \e
\]
for all $t \in \R$.
We wish to thank Quanhua Xu, who pointed out this trick to the second author.

Since  ($\overline M_{\p}, m_{\check  \p}$) is a compatible pair of
finite-rank maps  on  $(C^*_{red}(G),  L_1(VN(G)))$, this induces a 
well-defined
finite-rank map $T$ on
$C^*_{red}(G) +  L_1(VN(G))$ such that $(T\otimes id_\infty) \circ f 
\in \F$ and
$T\otimes id_\infty (a) = T\otimes id_\infty(f ({\frac 1p}))$.
Then we have
\begin{eqnarray*}
\|T\otimes id_\infty (a) &-& a \|_{L_p(VN(G)) \ten K_\infty}
\le   \|(T\otimes id_\infty)\circ f - f\|_\F \\
&=& \max \{\sup_{t\in \R}\{\|\overline M_{\p} \otimes id_\infty(f(it))
- f(it)\|_{C^*_{red}(G) \ten K_\infty}\}, \\
&{~} ~&
~ \, ~ \, ~ ~ \, ~ \, ~ ~ \, ~ \, ~
\sup_{t\in \R}\{\|m_{\check  \p} \otimes id_\infty(f(1+it)) -
f(1+it)\|_{ L_1(VN(G))  \ten K_\infty}\} \}
\le \e.
\end{eqnarray*}
This completes the proof.  \endproof

Using a similar argument, we can easily prove the following proposition.

\begin{proposition}
\label {P3.cb}
Let $G$ be a weakly amenable discrete group and $1 < p < \infty$.
Then $L_p(VN(G))$ has the CBAP with
\[
\Lambda(L_p(VN(G))) \le \Lambda(G).
\]
\end{proposition}
\begin{proof}
If $G$ is weakly amenable,
then there exists a net of $\{\p_\al\}$ in $A_c(G)$ such that
\[
\|m_{{\check  \p}_\al}\|_{cb} = \|\overline {M}_{\p_\al}\|_{cb} \le \Lambda(G)
\]
and  $m_{{\check  \p}_\al} \to id_{ L_1(VN(G))}$ and
$\overline M_{\p_\al} \to id_{C^*_{red}(G)}$
in the point-norm  topology.
For each $\al$, it is known from the complex interpolation theory
(see Pisier \cite [Proposition 2.1]{PiOH}) that
$(\overline M_{\p_\al}, m_{{\check  \p}_\al})$ induces a completely bounded
finite-rank  map $T_\alpha$ on $L_p(VN(G))$ such that
\[
\|T_\alpha\|_{cb} \le \|\overline M_{\p_\al}\|_{cb}^{1- \frac 1p}
\|m_{{\check  \p}_\al}\|_{cb}^{\frac 1p} \le \Lambda(G).
\]
The fact that   $T_\al \to id_{L_p(VN(G))}$ in the
point-norm topology follows from a similar argument that given
in the proof of Theorem 1.1.
\end{proof}

\section{Completely Contractive Approximation Property for $L_p(VN(G))$}

Given an operator space $V$, we let $\F(V, V)$ denote the space of all
(completely) bounded finite-rank maps on $V$.
Then we may identify the algebraic tensor product $V^* \otimes V$
with $\F(V, V)$  via the map
\[
u \in V^* \otimes V \to F_u = \sum f_i \otimes x_i.
\]
The \emph{$\lambda$-ball} of
$\F(V, V)$ is refered to as the set of all elements in $T\in \F(V, V)$ with
$\|T\|_{cb}\le \lambda$.
We first need  the following  lemma.

\begin{lemma}
\label {P4.1}
Let $V$ be an operator space and $1\le \la < \infty$.
Then $\Lambda (V) \le \la$ if and only if
for every $\e > 0$, $id_V$ is contained in the stable point-norm closure of
the  $(1+\e) \la $-ball of $\F(V, V)$.

\end{lemma}
\begin{proof}
$\Leftarrow :$
Given any finite subset  $S=\{v_1, \cdots, v_n\}$ of $V$ and $\e > 0$,
we let $I = \{\al\}$ denote
the index set consisting of all such $\al= \{S, \e\}$.
Then there is a canonical partial order on $I$ given by
$\al \preceq \al '$ if and only if $\e' \le \e$ and $S \subseteq S'$.

Given any $\al = \{S, \e\}\in I$, we let
\[
M = {\rm max} \{\|v_i\|: ~ v_i \in S\} +1.
\]
Since $id _V$ is contained in the stable point-norm closure of the
$(1+\frac \e M)\lambda $-ball of $\F(V, V)$,   there exists a finite-rank map
$T : V \to V$ such that
$\|T\|_{cb} \le (1+ \frac {\e} M) \lambda$ and
$\|T(v_i) - v_i\| < \e$ for all $v_i \in S$.
Then $T_\al ={ \frac 1 {1+ \frac {\e} M } } T  $ is a
finite-rank map on $V$ such that  $\|T_\al\|_{cb} \le \lambda$ and
\[
\|T_\al (v_i) - v_i\| \le {\frac {\frac \e M} {1+\frac \e M}}\|T(v_i)\| +
\|T(v_i) - v_i\| < (\la +1) \e.
\]
This shows that  $T_\al \to id_V$ in the point-norm topology, and thus
$\Lambda(V) \le \lambda$.

$\Rightarrow :$ This is obvious since if $\Lambda(V) \le \la$,
then there exists a net of finite-rank maps $T_\alpha$ on $V$ such that
$\|T_\al\|_{cb} \le \lambda$ and $T_\al \to id_V$ in the point-norm
topology.  This implies that $T_\al\to id_V $ in the stable point-norm
topology, and thus $id_V$ is contained in the stable point-norm
closure of the  $\lambda$-ball of $\F(V, V)$.
\end{proof}

Given operator spaces $V$ and $W$, we have the (complete) isometric inclusions
\[
\CB (V, W) \hookrightarrow \CB(V, W^{**}) \cong (V\hat \otimes W^*)^*.
\]
Then every $u \in V \hat \otimes W^*$ determines a bounded linear functional
$F_u$ on $\CB (V, W)$ with
\[
\|F_u\| \le \|u\|_{V \hat \otimes W^*}.
\]
To see this, suppose that we are given an  element $u \in V\hat \otimes W^*$.
We  define $F_u : \CB(V, W) \to {\Bbb C}$ to be given by
\[
F_u(T) = \al [g_{kl}(T(v_{ij}))]\be
\]
if  we can write
$u  = \al (v \otimes g)\be$ for some $\al \in M_{1, \infty^2}$,
$v = [v_{ij}] \in V  \ten K_\infty$, $g = [g_{kl}]\in W^* \ten K_\infty$ and
$\be \in M_{\infty^2, 1}$.
It is clear that  this is a well-defined linear functional on $\CB (V, W)$ and
\[
|F_u(T)| \le \|\al\| \, \|v\|\, \|g\|\, \|\be \| \, \|T\|_{cb}
\]
for all $T \in \CB (V, W)$.
This  shows that
\[
\|F_u\| \le \inf \{\|\al\|\, \|v\|\, \|g\|\, \|\be\|: ~
u = \al (v \otimes g)\be  \} = \|u\|_{V \hat \otimes W^*}.
\]
Moreover, it was shown  in \cite {ERmapping} that
these maps $F_u$ with $u \in V \hat \otimes W^*$ are exactly the
continuous linear functionals  with respect to the
stable point-norm topology on $\CB (V, W)$.

\begin{theorem}
\label {P4.2}
Let $R$ be a von Neumann algebra with a normal faithful tracial state
and with the QWEP, and let $1 < p < \infty$.
If $V$ is a complemented operator subspace
of $L_p(R)$  and $V$ has the OAP, then $V$ has the CBAP such that
\[
\Lambda (V) \le \|P\|_{cb},
\]
where $P$ is  a completely bounded projection from  $L_p(R)$ onto $V$.
\end{theorem}
\begin{proof}
 From Lemma \ref {P4.1}, it suffices to show that for every $\e >0$
the identity map $id_V$ is  contained in the stable point-norm
closure of the  $(1+\e) \|P\|_{cb}$-ball of  $\F(V, V)$.
Suppose not, i.e. suppose that there exists an $\e > 0$ such that
$id_V$ is not contained in the stable point-norm closure of the
$(1+{\e })\|P\|_{cb}$-ball of  $\F(V, V)$.
Then there exists a stable point-norm continuous linear functional $F$ on
$\CB (V, V)$ such that $F (id_V) = 1$ and $|F(T)| \le 1$ for all
$T \in \F(V, V)$ with $\|T\|_{cb} \le (1+{\e})\|P\|_{cb}$.
The latter condition is  equivalent to saying that we have
\[
|F(T)| \le \frac {\|T\|_{cb}} {(1+{\e})\|P\|_{cb}}
\]
for all $T \in \F(V, V)$.
Since $V^*\ten V$ can be identified with the cb-norm closure of
$V^* \otimes V\cong \F(V, V)$ in $\CB(V, V)$,
$F$ can be identified with a bounded linear functional on $V^*\ten V$,
which is still denoted by $F$,
with  norm
\[
\|F\|_{V^*\ten V}  \le \frac 1 {(1+{\e})\|P\|_{cb}} .
\]

We let ${\mathcal I} (V^*, V^*)$ denote the space of all
\emph{ completely integral maps} from $V^*$ into $V^*$,
and let ${\mathcal N} (V^*, V^*)$ denote the space of all
\emph {completely nuclear maps} from $V^*$ into $V^*$
(see details in \cite{ERmapping} and \cite {ERbook}).
For  $1 < p < \infty$, $L_p(R)$ is a reflexive space.
Then the closed subspace $V$ of  $L_p(R)$ is also reflexive, and thus
we have the isometry
\[
{\mathcal I} (V^*, V^*)  = (V^*\ten V)^*.
\]
Since $V$ has the OAP,  we have the isometry
\[
{\mathcal N} (V^*, V^*) = V^{**}\hat \otimes V^* =  V\hat \otimes V^*.
\]
Together with the QWEP condition on $R$,
Junge \cite{Ju} proved that  we actually  have the isometry
\[
{\mathcal I}(V^*, V^*)= {\mathcal N}(V^*, V^*) = V \hat \otimes V^* .
\]
Then we may choose an element  $u \in V \hat \otimes V^*$ such that
$F = F_u $ and
\[
\|u\|_{V \hat \otimes V^*} =\|F\| \le
  \frac 1 { (1+\e)\|P\|_{cb}}.
\]

Since $V$ has the OAP, there exists a net of finite-rank maps
$T_\al \in \F(V, V)$ such that $T_\al \to id_V$ in the stable point-norm
topology.
Then we obtain a contradition since
\begin{eqnarray*}
1 = F(id_V) &=& \lim_\al |F (T_\al)| = \lim_\al |F_u( T_\al)|\\
&=  & |F_u(id_V)| \le \|u\|_{V\hat \otimes V^*} \|id_V\|_{cb}
\le \frac 1 {(1+\e)\|P\|_{cb}} < 1.
\end{eqnarray*}
This   shows that $\Lambda(V) \le \|P\|_{cb}$.
\end{proof}

Since CBAP implies OAP, we may obtain the following proposition
from Theorem \ref {P4.2} by considering $V = L_p(R)$ and $P = id_{L_p(R)}$.

\begin{proposition}
\label {P4.2b}
Let $R$ be a von Neumann algebra with a normal faithful tracial state
and with the QWEP.
For $1 < p < \infty$, the following are equivalent:

\begin{itemize}
\item [(1)] $L_p(R)$ has the CCAP,

\item [(2)] $L_p(R)$ has the CBAP,

\item [(3)] $L_p(R)$ has the OAP.
\end{itemize}
\end{proposition}

\noindent \emph{Proof of Theorem 1.2.}
Let $G$ be a discrete group with the AP and let $VN(G)$ have the QWEP.
Then for $1 < p < \infty$, $L_p(NV(G))$ has the OAP by Theorem 1.1,
and thus has the CCAP by Proposition \ref {P4.2b}.
$\square$

Let $G$ be a discrete group.
Then there is a  unital normal $*$-isomorphic
injection $\pi: VN(G) \to VN(G \times G)$, which is defined by
\begin{equation}
\label{F4.pi}
\pi (\lambda(s)) = \lambda(s) \otimes \lambda(s),
\end{equation}
and is known as the \emph{co-multiplication} of
the \emph{Hopf von Neumann algebra}  $VN(G)$.
The map $\pi$ is  invariant with respect to the tracial states on
$VN(G)$ and $VN(G\times G)$, i.e. we have
\begin{equation}
(\tau \otimes \tau)\circ \pi (x) =  \tau (x)
\end{equation}
for all $x \in VN(G)$.
We may use $\pi$ to identify $VN(G)$ with the von Neumann
subalgebra $\pi(VN(G))$ of $VN(G\times G)$.
Since $\tau\otimes \tau$, restricted to $\pi(VN(G))$,
is again a normal faithful tracial state on $\pi(VN(G))$,
there is a unique $\tau \otimes \tau$-invariant normal conditional
expectation  ${\mathcal E}$ from  $VN(G\times G) $ onto  $\pi(VN(G))$
(see \cite {St}) which satisfies
\[
{\mathcal E}(\sum_{s, t} \al_{s, t}\la(s) \otimes \la(t)) =
\sum_{s} \al_{s, s}\la(s) \otimes \la(s)
\]
for all $\sum_{s, t} \al_{s, t} \la(s) \otimes \la(t) \in VN(G\times G)$.
If  we let
\begin{equation}
\rho = \pi^{-1} \circ {\mathcal E} :  VN(G \times G) \to VN(G),
\end{equation}
then $\rho$ is a normal complete contraction such that
$\rho \circ \pi = id_{VN(G)}$,
and
\[
\tau \circ \rho (y) = \tau \otimes \tau ({\mathcal E}(y))
= \tau \otimes \tau (y)
\]
for all $y\in VN(G \times G)$.
We note that
\begin{equation}
\label{F4.rho}
\kappa \circ \rho \circ (\kappa \otimes \kappa) = \rho
\end{equation}
on  $VN(G\times G)$ since for every
$x = \sum_{s, t} \al_{s, t} \la(s) \otimes \la(t) \in VN(G\times G)$ we have
\begin{eqnarray*}
\kappa \circ \rho \circ (\kappa \otimes \kappa) (x) &= &
\kappa \circ \rho (\sum_{s, t} \al_{s, t} \la(s^{-1}) \otimes \la(t^{-1}))\\
&=& \kappa (\sum_{s} \al_{s, s} \la(s^{-1}))\\
&=& \sum_{s} \al_{s, s} \la(s) = \rho(x) .
\end{eqnarray*}
The pre-adjoint map $\rho_*: A(G) \to A(G\times G)$ is a completely
isometric inclusion, and induces a complete isometry
\[
\check \rho_* = (\kappa_*\otimes \kappa_*)\circ \rho_* \circ \kappa_*:
L_1(VN(G)) \to L_1(VN(G\times G)).
\]
It follows from (\ref {F4.rho}) that we may identify $\check \rho_*$ with
$\rho_*$ and  obtain
\begin{eqnarray*}
\langle\check \rho_* (j_1(\la(s)), ~ x\rangle
&=& \langle j_1(\la(s)), ~  \rho (x) \rangle
= \tau( \la(s) \pi^{-1}\circ{\mathcal E} (x)) \\
&=& \tau \otimes \tau (\pi(\la(s)) x)
=\langle  (j_1 \otimes j_1) \circ \pi (\la(s)), ~ x \rangle
\end{eqnarray*}
for all $s\in G$ and $x \in VN(G\times G)$.
This shows that if we let
\[
\pi_\infty = \pi_{|C^*_{red}(G)}: C^*_{red}(G)\to C^*_{red}(G\times G)
\]
and
\[
\pi_1 = \check \rho_* :  L_1(VN(G))  \to L_1(VN(G\times G)) ,
\]
then
$(\pi_\infty, \pi_1)$ satisfies
\[
\pi_1\circ j_1 = (j_1 \otimes j_1) \circ \pi_\infty,
\]
and thus is a compatible pair of
complete contractions from
the compatible couple $(C^*_{red}(G), L_1(VN(G)))$ into the compatible couple
$(C^*_{red}(G\times G), L_1(VN(G\times G)))$.
By complex interpolation,
we obtain,  for each $1< p < \infty$, a
complete contraction
\[
\pi_p : L_p(VN(G)) \to L_p(VN(G\times G)),
\]
which satisfies
\begin{equation}
\pi_p(j_p(\la(s)) = j_p(\lambda(s)) \otimes j_p(\lambda(s)).
\end{equation}

On the other hand, we let
\[
\rho_\infty = \rho_{|C^*_{red}(G\times G)}: C^*_{red}(G\times G) \to 
C^*_{red}(G)
\]
and
\[
\rho_1= \check \pi_* : L_1(VN(G\times G)) \to  L_1(VN(G)),
\]
where $\check \pi_* = \kappa_* \circ \pi_* \circ (\kappa_* \otimes \kappa_*)$.
Then we may identify $\check \pi_* $ with $\pi_*$
and show that $(\rho_\infty, \rho_1)$ is a
compatible pair of complete contractions from
the compatible couple  $(C^*_{red}(G\times G), L_1(VN(G\times G)))$ into the
compatible couple $(C^*_{red}(G), L_1(VN(G)))$.
By complex interpolation, we obtain, for each $1 < p < \infty$, a complete
  contraction
\[
\rho_p: L_p(VN(G\times G)) \to L_p(VN(G)),
\]
which satisfies that $\rho_p \circ \pi_p = id_{L_p(VN(G))}$.
This shows that $\pi_p$ is actually a completely isometric embedding from
$L_p(VN(G))$ into $ L_p(VN(G\times G))$, and the range space
$\pi_p(L_p(VN(G))$ is completely contractively
complemented in  $L_p(VN(G\times G))$.

Let $T$ be a  completely bounded map on $L_p(VN(G))$.
If $VN(G)$ has the QWEP, Junge \cite {Ju} proved
  that $T$ determines a completely bounded map
$id_{ L_p(VN(G))}\otimes T$  on $L_p(VN(G\times G))$, which
satisfies
\[
\|id_{ L_p(VN(G))}\otimes T\|_{cb} \le \|T\|_{cb}.
\]
 From this, we may construct another completely bounded map
\begin{equation}
\label{F4.HK}
\tilde T = \pi_{p'}^* \circ (id_{ L_p(VN(G))}\otimes T)\circ \pi_p
\end{equation}
on $L_p(VN(G))$.
We note that this  construction  is motivated by an argument given in
\cite [Theorem 2.1]{HK}.
Using Theorem 1.2 and a similar technique used in \cite [Theorem 2.1]{HK},
we may obtain the  following result.

\begin{theorem}
\label {P4.multiplier}
Let  $1 < p < \infty$.
If  $G$ is a discrete group with the AP and $VN(G)$ has the QWEP,
then there exists a net of completely
contractive finite-rank maps  $T_\al : L_p(VN(G)) \to L_p(VN(G))$
such that  $T_\al \to  id_{L_p(VN(G))}$ in the
point-norm topology, and for each $\al$,
\[
T_\al (j_p(\lambda(s))) = \p_\al(s) j_p(\lambda(s))
\]
for some  $\p_\alpha$ in $A_c(G)$.
\end{theorem}
\begin{proof}
With the hypothesis, we know from Theorem 1.2 that $L_p(VN(G))$ has the
CCAP.
Let $\{T_\alpha\}$ be a net of
finite-rank complete contractions on $L_p(VN(G))$ such that
$T_\al \to id_{L_p(VN(G))}$ in the point-norm topology.
Since $j_p(\la(\C [G]))$ is norm dense in $L_p(VN(G))$,
we may assume, without loss of generality,  that the range of
$T_\al$ is contained in $j_p(\la(\C [G]))$
(consider a  small perturbation if necessary).
Using (\ref {F4.HK}), we get a net of completely bounded maps
\[
\tilde T_\al =  \pi_{p'}^* \circ (id_{ L_p(VN(G))}\otimes T_\al)\circ \pi_p
\]
on $L_p(VN(G))$ such that  $\|\tilde T_\al\|_{cb} \le \|T_\al\|_{cb} \le 1$.

For each $\al$, we let  $\p_\al$ be the function on $G$ defined by
\[
\p_\al (s) = \langle T_\al (j_p(\lambda(s)),~  j_{p'}(\lambda (s^{-1}))\rangle.
\]
Since the range of $T_\al$ is contained in $j_p(\la(\C [G]))$,
we can write
\[
T_\al = \sum _i f^\al_i \otimes j_p(\la(s^\al_i))
\]
for some $f^\al_i \in L_p(VN(G))^*$ and $s^\al_i \in G$.
Then it is easy to see that
$\p_\al $ has a finite support and thus is an element in $A_c(G)$.
Moreover,  we have
\[
\tilde T_\al (j_p(\lambda(s))) = \p_\al(s)j_p(\lambda(s))
\]
because
\begin{eqnarray*}
\langle  \tilde T_\al(j_p(\lambda(s))),~  j_{p'}(\lambda(t))\rangle
&=& \langle  (id_{ L_p(VN(G))}\otimes T_\al)\circ \pi_p(j_p(\lambda(s))), ~
\pi_{p'}(j_{p'}(\lambda(t))) \rangle  \\
&=&  \langle j_p(\lambda(s)) \otimes T_\al(j_p(\lambda(s))), ~
j_{p'}(\lambda(t)) \otimes j_{p'}(\lambda(t)) \rangle  \\
&=& \langle T_\al(j_p(\lambda(s)), ~ j_{p'}(\lambda (t) \rangle\,
\tau(\lambda(s) \lambda(t) ) \\
&=& \langle T_\al(j_p(\lambda(s)), ~ j_{p'}(\lambda (s^{-1}) \rangle\,
\tau(\lambda(s) \lambda(t) ) \\
&=&   \langle \p_\al(s) j_p(\lambda(s)), ~ j_{p'}(\lambda(t)) \rangle
\end{eqnarray*}
hold for all  $s, t\in G$.
This shows that for each $\al$,
$\tilde T_\al = m_{\p_\al}$ is a  finite-rank multiplier map
on $L_p(VN(G))$.
Since $T_\al \to id_{L_p(VN(G))}$ in the point-norm topology, we have
\[
\p_\al(s) = \langle T_\al (j_p(\lambda(s)),~  j_{p'}(\lambda (s^{-1}))\rangle
\to \langle j_p(\lambda(s)),~  j_{p'}(\lambda (s^{-1}))\rangle
= \tau(\la(s)\la(s^{-1}))= 1
\]
for all $s \in G$.
This implies that
$\tilde T_\al \to id_{L_p(VN(G))}$ in the point-norm topology.
\end{proof}

\section{Local Structure and Completely Bounded Schauder Basis}

In a recent paper \cite {JNRX}, we studied the complemented
${\mathcal O \mathcal L}_p$ spaces and completely bounded Schauder
basis for operator spaces.
We recall that an operator space $V$ is said to be a
${\mathcal C \mathcal O \mathcal L}_p$ space if there exists a
constant $\lambda \ge 1$ such that for every finite-dimensional
subspace $E \subseteq V$ there exists a finite-dimensional $C^*$-algebra
$\A$ and a diagram of completely bounded maps
\[
\begin{array}{ccccc}
&  & V &  &  \\
& {\scriptstyle r}\nearrow  &  & \searrow {\scriptstyle s} &  \\
L_p(\A) &  & \stackrel{id }{\longrightarrow } &  & L_p(\A)
\end{array}
\]
with $\|r\|_{cb} \|s\|_{cb} \le \lambda$ and $E \subseteq r(L_p(\A))$.
This definition is slightly stronger than the notion of
${\mathcal O \mathcal L}_p$ spaces introduced  in \cite {ER2}.

A sequence of elements $\{x_n\}$ in a separable operator space $V$ is a
\emph{ Schauder basis}  for $V$ if every element $x \in V $ can be
uniquely written as
$x  = \sum_{n = 1} \alpha_n x_n$,
where the sum converges in the norm topology.
It is well-known from Banach space theory
that if $\{x_n\}$ is a basis for $V$, then
the natural projections
\[
P_k : x \in V \to \sum_{n=1}^k \alpha_n x_n ~ \in ~
\text{span}\{x_1, \cdots, x_k\}
\]
are uniformly bounded, i.e. $\sup\{\|P_k\|\} < \infty$.
We say that $\{x_n\}$ is a \emph{completely bounded Schauder basis} for $V$ if
$\sup\{\|P_k\|_{cb}\} < \infty$.
It was shown in \cite {JNRX} that if
$R$ is an injective von Neumann algebra with a separable predual, then
for $1\le p < \infty$,
$L_p(R)$ is a  ${\mathcal C \mathcal O \mathcal L}_p$ space and
has a completely bounded Schauder basis.
In general, if $R$ is a von Neumann
algebra with a separable predual and has the QWEP,
it was shown in  \cite {JNRX} for $1 < p < \infty$,
$L_p(R)$ has the CBAP  if and only if $L_p(R)$
is a  ${\mathcal C \mathcal O \mathcal L}_p$ space
(respectively, $L_p(R)$ has a completely bounded
Schauder basis).
As a consequence of this result and Theorem 1.2 (or Proposition \ref {P4.2b}),
we can conclude  the following  result for discrete groups.

\begin{theorem}
\label{P4.T}
Let $G$ be a countable discrete group and $1 < p < \infty$.
If $G$ has  the AP and $VN(G)$ has the QWEP,
then $L_p(VN(G))$ is a  ${\mathcal C \mathcal O \mathcal L}_p$ space
and  has a completely bounded  Schauder basis.
\end{theorem}

It is known that for $2 \le n < \infty$, the free  group of $n$ generators
${\Bbb F}_n$ is weakly amenable with $\Lambda({\Bbb F}_n) = 1$
(see \cite {DH} and  \cite {Ha1}) and
$VN({\Bbb F}_n)$ has the QWEP  (see \cite {Ki1}).
Then $L_p(VN({\Bbb F}_n))$ is a ${\mathcal C \mathcal O \mathcal L}_p$ space
and  has a completely bounded  Schauder basis.
In general, we wish to consider residually finite discrete groups with the AP.
We recall that a discrete group $G$ is \emph{residually finite} if
for any finitely many distinct elements $s_1, \cdots, s_n \in  G$, there is a
homomorphism $\theta$ from $G$ onto
a finite group $H$ such that
the elements $\theta(s_1), \cdots, \theta(s_n)$ are distinct in $H$.
It was shown by Wassermann \cite {Wa} that every residually finite
discrete group $G$ has the \emph {property} ($F$) (see definition
in \cite {Wa} and \cite {Ki1}), and thus its group von Neumann algebra
$VN(G)$ has the QWEP (see Kirchberg \cite [$\S 7$]{Ki1}).
 From this, we can conclude the following proposition.

\begin{proposition}
\label {P4.RF}
Let $G$ be a countable residually finite discrete group with the AP, then
for $1 < p < \infty$,  $L_p(VN(G))$ has the CCAP, and thus is a
${\mathcal C \mathcal O \mathcal L}_p$ space
and  has a completely bounded  Schauder basis.
\end{proposition}

In the following let us discuss some examples of countable residually finite
discrete groups for which the non-commutative $L_p(VN(G))$ spaces are
${\mathcal C \mathcal O \mathcal L}_p$ spaces
and  have  completely bounded  Schauder bases.

\begin{example}
{\rm
Since the integer valued special linear group
$SL(2, {\Bbb Z})$ is a closed discrete subgroup of $SL(2, {\Bbb R})$,
it is known from \cite {DH} that  $SL(2, {\Bbb Z})$ is weakly amenable
with $\Lambda(SL(2, {\Bbb Z})) = 1$.
This group is also residually finite since
for any finitely many distinct elements $s_1, \cdots, s_n \in
SL(2,  {\Bbb Z})$, there is a sufficiently   large prime number $p$
and a   natural homomorphism $\theta_p$
from $SL(2,  {\Bbb Z})$ onto
the finite group $SL(2,  {\Bbb Z}_p)$,
where ${\Bbb Z}_p = {\Bbb Z}/p{\Bbb Z}$,
such that the elements $\theta_p(s_1), \cdots, \theta_p(s_n)$
are distinct in $SL(2,  {\Bbb Z}_p)$}.
\end{example}

\begin{example}
{\rm Let us consider the discrete group
${\Bbb Z}^2 \rtimes SL(2,  {\Bbb Z})$, the semi-direct
product of ${\Bbb Z}^2$ by $SL(2,  {\Bbb Z})$.
It is known from \cite {HK} that this group
has the AP, but is not weakly amenable.
On the other hand, it is easy to see that this group is residually finite
since we may  consider  the natural homomorphisms $\theta_p$
from $ {\Bbb Z}^2 \rtimes SL(2,  {\Bbb Z})$ onto
the finite groups ${\Bbb Z}_p^2 \rtimes SL(2,  {\Bbb Z}_p)$
for prime numbers $p$}.
\end{example}

\begin{example}
{\rm  It is known (by Borel and Harish-Chandra \cite {BC}) that
every arithmetic subgroup $\Gamma$ of $Sp(1, n)$ (with $n \in {\Bbb N}$) is
a \emph{lattice} of $Sp(1, n)$, i.e. a closed discrete subgroup of $Sp(1, n)$
for which $Sp(1, n)/\Gamma$ has a bounded $Sp(1, n)$-invariant measure.
In this case, $\Gamma$ is weakly amenable with
\[
\Lambda(\Gamma) = \Lambda(Sp(1, n)) = 2n-1
\]
(see \cite {Ha2} and \cite {CH}).
Let  ${\Bbb H}_{int}={\Bbb Z}+{\Bbb Z} i+{\Bbb Z}j+{\Bbb Z}k$
denote the quaternionic integers ring, and let $\Gamma_n$ denote the
subgroup of $Sp(1, n)$ consisting of all
$(n+1)\times (n+1)$ matrices with values in ${\Bbb H}_{int}$,
which leave  the sesquilinear form
  \[
Q(x,y) =  \bar{y}_0x_0-\sum\limits_{m=1}^n \bar y_m {x}_m ~\, ~\,  ~\, ~\,
~\, ~\, ~\, ~\,~\, ~\,
(x_m, y_m \in {\Bbb H}_{int})
\]
invariant.
Then $\Gamma_n$ is a lattice of $Sp(1, n)$ such that
$\Lambda(\Gamma_n) = \Lambda(Sp(1, n)) = 2n-1$ (see \cite {CH} $\S 6$).
Considering the canonical  representations of $\Gamma_n$ into the finite ring
$M_{n+1}({\Bbb Z}_p + {\Bbb Z}_p  i+{\Bbb Z}_p i+ {\Bbb Z}_p k)$
for prime numbers $p$, we see that $\Gamma_n$ is residually finite,
and thus $VN(\Gamma_n)$  has the QWEP}.
\end{example}

\begin{example}
{\rm Let $\Gamma_n$ be the lattice of  $Sp(1, n)$ discussed in (3),
and let $G =\amalg_n \Gamma_n$ be the discrete subgroup of
the infinite product group
$\prod_n \Gamma_n$, which  consists of all sequences $g=(g_n)$ with
$g_n \in \Gamma _n$ such that  all but finitely many $g_n$ are equal to
the unit element $1_n$ of $\Gamma_n$.
Then for each $n \in {\Bbb N}$, we may identify the normal subgroup
\[
G_n = \{(g_i): ~ g_i = 1_i ~ {\rm for}  ~ i > n\}
\]
of $G$ with the finite product group
$\Gamma_1\times  \cdots \times \Gamma_n$, and
identify
\[
VN(G_n)= VN(\Gamma_1) \overline \otimes \cdots \overline \otimes VN(\Gamma_n)
\]
with a von Neumann subalgebra of $VN(G)$.
Since $\{G_n\}$ is an increasing sequence of normal subgroups in $G$ and
$G = \bigcup_n G_n$,
we get an increasing sequence of finite von Neumann algebras
\[
VN(G_1) \hookrightarrow \cdots \hookrightarrow VN(G_n) \hookrightarrow
\cdots \hookrightarrow
VN(G),
\]
which preserve the tracial states and have a weak$^*$ dense
union  $\bigcup_n VN(G_n)$  in $VN(G)$.
Then
there is an increasing  sequence of (completely contractive)
normal conditional expectations ${\mathcal E}_n$ from  $VN(G)$  onto $VN(G_n)$
such that  ${\mathcal E}_n \to id_{VN(G)}$ in point-weak$^*$ topology (and thus
in stable point-weak$^*$ topology).

It was shown by Cowling and Haagerup \cite {CH} that the groups  $G_n$
are weakly  amenable with
\[
\Lambda (G_n) = \Lambda(VN(G_n))= \Lambda(VN(\Gamma_1))\cdots
\Lambda(VN(\Gamma_n)) = 1 \cdot 3 \cdots (2n-1).
\]
Then we can conclude that
\[
\Lambda(G)= \sup \{\Lambda(G_n)\} = \infty.
\]
This shows that the group $G$ is not weakly amenable.
But $G$ has the AP.
This follows from the fact that any point $u \in VN(G) \overline 
\otimes M_\infty$
can approximated  by $\{({\mathcal E}_n\otimes id_\infty)(u)\}$ in
the weak$^*$ topology on
$VN(G) \overline \otimes M_\infty$,
and for each $n$, $VN(G_n)$ has the weak$^*$ CBAP (and thus weak$^*$ OAP).
Moreover, for each $n\in {\Bbb N}$,
$G_n \cong \Gamma_1\times  \cdots \times \Gamma_n$ is residually 
finite and thus
$G$ is residually finite.
This provides another example of a discrete residually
finite group which has the AP, but is not weakly amenable. }
\end{example}

\end{document}